\documentclass{scrartcl}
\usepackage{typearea}
\typearea{12}
\usepackage{amssymb, amsmath, amsthm}
\usepackage{enumitem}
\usepackage[all, warning]{onlyamsmath}
\usepackage{graphicx}
\usepackage[ruled, lined]{algorithm2e}
\usepackage{tikz}\usetikzlibrary{arrows}
\usepackage{color}
\usepackage{ulem}
\usepackage{authblk}
\usepackage{url}
\usepackage{pdflscape} 

\newtheorem{THEO}{Theorem}[section]
\newtheorem{ALGo}[THEO]{Algorithm}
\newtheorem{ASSM}[THEO]{Assumption}
\newtheorem{CONJ}[THEO]{Conjecture}
\newtheorem{COND}[THEO]{Condition}
\newtheorem{CORO}[THEO]{Corollary}
\newtheorem{DEFI}[THEO]{Definition}
\newtheorem{EXAMP}[THEO]{Example}
\newtheorem{FACT}[THEO]{Fact}
\newtheorem{HYPO}[THEO]{Hypothesis}
\newtheorem{LEMM}[THEO]{Lemma}
\newtheorem{PROB}[THEO]{Problem}
\newtheorem{PROP}[THEO]{Proposition}
\newtheorem{REMA}[THEO]{Remark}
\newcommand{\theo}{\begin{THEO}}
\newcommand{\algo}{\begin{ALGo}}
\newcommand{\assm}{\begin{ASSM}}
\newcommand{\cond}{\begin{COND}}
\newcommand{\conj}{\begin{CONJ}}
\newcommand{\coro}{\begin{CORO}}
\newcommand{\defi}{\begin{DEFI}}
\newcommand{\examp}{\begin{EXAMP}}
\newcommand{\fact}{\begin{FACT}}
\newcommand{\hypo}{\begin{HYPO}}
\newcommand{\lemm}{\begin{LEMM}}
\newcommand{\prob}{\begin{PROB}}
\newcommand{\prop}{\begin{PROP}}
\newcommand{\rema}{\begin{REMA}}
\newcommand{\etheo}{\end{THEO}}
\newcommand{\ealgo}{\end{ALGo}}
\newcommand{\eassm}{\end{ASSM}}
\newcommand{\econd}{\end{COND}}
\newcommand{\econj}{\end{CONJ}}
\newcommand{\ecoro}{\end{CORO}}
\newcommand{\edefi}{\end{DEFI}}
\newcommand{\eexamp}{\end{EXAMP}}
\newcommand{\efact}{\end{FACT}}
\newcommand{\ehypo}{\end{HYPO}}
\newcommand{\elemm}{\end{LEMM}}
\newcommand{\eprob}{\end{PROB}}
\newcommand{\eprop}{\end{PROP}}
\newcommand{\erema}{\end{REMA}}
\DeclareMathOperator{\He}{He}
\DeclareMathOperator{\rank}{rank}

\makeatletter
\global\let\tikz@ensure@dollar@catcode=\relax
\makeatother

\title{Application of Facial Reduction to \\ $H_\infty$ State Feedback Control Problem}
\author[1]{Hayato Waki\thanks{744 Motooka, Nishi-ku, Fukuoka 819-0395, Japan. waki@imi.kyushu-u.ac.jp}}
\author[2]{Noboru Sebe\thanks{680-4 Kawazu, Iizuka-shi, Fukuoka 820-8502, Japan.  sebe@ai.kyutech.ac.jp}}
\affil[1]{Institute of Mathematics for Industry, Kyushu University}
\affil[2]{Department of Artificial Intelligence, 
Faculty of Computer Science and Systems Engineering, 
Kyushu Institute of Technology}
\date{\today}

\begin{document}
\maketitle

\begin{abstract}
One often encounters numerical difficulties in solving linear matrix inequality (LMI) problems obtained from $H_\infty$ control problems. We discuss the reason from the viewpoint of optimization, and 
provide necessary and sufficient conditions for LMI problem and its dual not to be strongly feasible. Moreover, we  interpret them in terms of control system. In this analysis, facial reduction, which was proposed by Borwein and Wolkowicz, plays an important role. We show that a necessary and sufficient condition  closely related to the existence of invariant zeros in the closed left-half plane in the system, and present a way to remove the numerical difficulty with the null vectors associated with invariant zeros in the closed left-half plane. 
Numerical results show that the numerical stability is improved by applying it. 

\end{abstract}

\section{Introduction}\label{sec:intro}

$H_\infty$ control problems have attracted attention from a lot of researchers in control and optimization fields since primal-dual interior-point methods (PDIPMs) were proposed in 90's.  $H_\infty$ control problems can be reformulated as linear matrix inequality (LMI) problem and be efficiently solved 
 by LMI software, such as SeDuMi \cite{Sturm99}, SDPT3 \cite{Toh99} and SDPA \cite{Fujisawa03}, etc. 
 Still, one often encounters some numerical difficulties in solving  LMI problems obtained from $H_{\infty}$ control problems by these LMI software. 
 
 The purpose of this manuscript is to investigate the reason why LMI software often return  inaccurate solutions and values to a fixed tolerance. It is empirically known that when either an LMI problem or its dual is not strongly feasible, the numerical instability will occur. In that case, optimal solutions may not exist. In addition, PDIPMs may not converge numerically. In fact, theoretical results of PDIPMs are required to be strongly feasible for both LMI problems and its dual.  
 
 Facial Reduction  (FR) is useful in such cases. It can detect whether a given convex optimization problem is strongly feasible or not. If not so, it finds a certificate that the problem is not strongly feasible, and generate an equivalent convex problem that is strongly feasible.   However, the execution of FR needs much more computation cost than solving the original LMI problem. 
 
 In this manuscript, we deal with state feedback controls for linear time invariant systems. 
 We present necessary and sufficient conditions for LMI problems and its dual of $H_\infty$ control problems obtained from them not to be strongly feasible. One of them is related to the existence of a stable invariant zero in a given system. In other words, the dual problem is not strongly feasible if the system has a stable invariant zero.  
 
We also provide how to remove  the numerical difficulty caused by non-strong feasibility. This is also based on FR and it generates a smaller LMI problem by using invariant zeros in the closed left-half plane. Interestingly, the resulting LMI problem can be obtained from a subsystem of the closed loop system.  In fact, a non-singular matrix used in the reduction consists of invariant zeros in the closed left-half plane and  plays an important role in FR for a given LMI problem. We show that the subsystem is obtained by applying the transformation with the matrix into the closed loop system. This implies that one can retrieve a state feedback gain for the closed loop from one for the subsystem. 
We also present numerical experiments to see the improvement on the numerical stability.
 
 \subsection{Literature related to this topic}
 \cite{Balakrishnan03} applies FR to the analysis and design of $H_{2}$ state feedback control, and showed the relationship between the strong feasibility and  invariant zeros in $H_{2}$ state feedback control. 
 \cite{Vandenberghe05} proposes a fast implementation for LMI problems based on KYP lemma. This approach reformulates the dual of a given  LMI problem into the form of the LMI problem.  Although this implementation reduces the computational complexity, this is essentially different from FR. In fact, this does not change the strong feasibility of the resulting LMI problem. 
   
\cite{Scherer92a, Scherer92b, Stoorvogel91} investigate the effect of invariant zeros for the performance index $\gamma$ in $H_{\infty}$ feedback control and provide the $H_{\infty}$ norm conditions without the assumptions on zeros of the systems. The provided conditions are the reduced size Riccati inequalities. We provide similar results by applying FR.
This implies that we should reduce LMI problems from the
view point of numerical accuracy, as far as using
software based on PDIPMs. Control systems whose $D_{12}$ is not full column rank has been also handled in this manuscript as well as those papers. We can verify by using FR that the optimal performance index $\gamma$ for such a system is equivalent to $\gamma$ for a system in which a differentiator is added.

The organization of this manuscript is as follows: In section \ref{sec:preliminary}, we provide some facts on LMI problems and FR. We interpret them with terms of systems and provide a reduction of LMI problems and closed loop systems in sections \ref{sec:main} and \ref{sec:reduction}. Section \ref{sec:experiment} provides a numerical experiment. We give a conclusion of this manuscript in section \ref{sec:conclusion}. This manuscript is based on \cite{Waki15}. We add some technical proofs for some results and a detailed numerical results in this manuscript. 

\subsection{Notation and symbols}
Let $\mathbb{R}$ and $\mathbb{C}$ be the sets of real and complex numbers, respectively. We represent the sets of complex numbers with nonnegative real parts and nonpositive real parts by $\overline{\mathbb{C}_+}$ and $\overline{\mathbb{C}_-}$, respectively. Let $\mathbb{R}^n$ be the set of $n$-dimensional Euclidean space. Let $\mathbb{R}^{m\times n}$, $\mathbb{S}^n$, $\mathbb{S}^{n}_+$  and $\mathbb{S}^n_{++}$ be the sets of $m\times n$ real matrices, $n\times n$ symmetric matrices, $n\times n$ positive semidefinite matrices and $n\times n$ positive definite matrices, respectively. We denote the $m\times n$ zero matrix by $O_{m\times n}$. For $A, B\in\mathbb{R}^{m\times n}$, we define $A\bullet B :=\mbox{Trace}(AB^T)$. For $A\in\mathbb{R}^{n\times n}$, we define $\He(A) = A+A^T$. 

\section{Preliminary}\label{sec:preliminary}
\subsection{Linear matrix inequality and its strong duality theorem}\label{subset:lmi}
Linear Matrix Inequality (LMI) problem is formulated as follows: 
\begin{align}
\label{LMI}
\theta_P^*&=\inf_{x\in\mathbb{R}^m, X\in\mathbb{S}^n} \left\{c^Tx : X = \sum_{j\in\mathcal{M}}x_j F_j - F_0, X\in\mathbb{S}^n_+\right\}, 
\end{align}
where $c\in\mathbb{R}^m$, $\mathcal{M}=\{1, \ldots, m\}$ and $F_0, \ldots, F_m\in\mathbb{S}^n$.  Throughout this manuscript, we assume that $F_1, \ldots, F_m$ is linearly independent. One can obtain an approximation of an optimal solution of (\ref{LMI}) to any given tolerance by applying PDIPMs. In fact, many variants of PDIPMs are proposed and implemented as optimization software in SeDuMi \cite{Sturm99}, SDPT3 \cite{Toh99}, SDPA \cite{Fujisawa03}, etc. 

The dual problem of (\ref{LMI}) can be formulated as follows: 
\begin{align}
\label{dualLMI}
\theta_D^*&=\sup_{Y\in\mathbb{S}^n} \left\{F_0\bullet Y : F_j\bullet Y = c_j \ (j\in\mathcal{M}), Y\in\mathbb{S}^n_+\right\}. 
\end{align} It is well-known that the strong duality theorem holds for (\ref{LMI}) and its dual (\ref{dualLMI}) under a mild assumption which is called {\itshape Slater's condition}. See Theorem \ref{sdt} below.  Unlike to Linear Program (LP), the strong duality theorem for (\ref{LMI}) and (\ref{dualLMI}) requires such a condition. Problem (\ref{LMI}) is said to be {\itshape strongly feasible} if there exists $(\hat{x}, \hat{X})\in\mathbb{R}^m\times\mathbb{S}^n_{++}$ such that $\hat{X}=\displaystyle\sum_{j\in\mathcal{M}} \hat{x}_jF_j-F_0$. Similarly, (\ref{dualLMI}) is said to be {\itshape strongly feasible} if there exists $\hat{Y}\in\mathbb{S}^n_{++}$ such that $F_j\bullet \hat{Y} = c_j$ for all $j\in\mathcal{M}$. Slater's condition holds in (\ref{LMI}) (resp., (\ref{dualLMI})) is satisfied if (\ref{LMI}) (resp., (\ref{dualLMI})) is strongly feasible.

\theo\label{sdt}(\cite[Theorem 3.2.8]{Renegar01}; see also \cite[Theorem 2.2]{deKlerk02} and \cite[Section 4.7]{Gartner12})
If (\ref{dualLMI}) is strongly feasible and (\ref{LMI}) is feasible, then $\theta_P^*=\theta_D^*$ and (\ref{LMI}) has an optimal solution. Similarly, if (\ref{LMI}) is strongly feasible and (\ref{dualLMI}) is feasible, then $\theta_P^*=\theta_D^*$ and (\ref{dualLMI}) has an optimal solution. 
\etheo
It should be noted that Slater's conditions for both (\ref{LMI}) and its dual (\ref{dualLMI}) guarantee the convergence of PDIPMs. See \cite{deKlerk02, Renegar01}  and references therein for more details. 

\subsection{Facial reduction for LMI problem (\ref{LMI}) and its dual (\ref{dualLMI})}\label{subset:fra}
Borwein and Wolkowicz \cite{Borwein81} propose an approach to the strong duality theorem for convex optimization problems without assuming Slater's condition and any constraint qualifications. The approach is called {\itshape Facial Reduction (FR)}. \cite{Ramana97, Ramana97SIOPT} discuss the application of FR into LMI problem (\ref{LMI}) and its dual (\ref{dualLMI}). Some extensions are discussed in 
\cite{Pataki13, Waki13}. 

Facial Reduction  (FR) is an algorithm that generates an LMI problem which is strongly feasible by using a given LMI problem, or that detects the infeasibility. FR has the property of a finite convergence.  In each iteration, FR finds a nonzero solution of a problem that consists of LMIs, or detects the infeasibility. Since one has to solve a similar LMI problem in each iteration, the computation spends as much cost as solving the original LMI problem. Moreover, FR requires an exact solution of a generated problem in each iteration, and thus FR is not practical algorithm from the viewpoint of computational practice. 

Still, when (\ref{LMI}) and/or its dual (\ref{dualLMI}) is not strongly feasible, one often encounters numerical difficulty on solving them. See \cite{Henrion05, Navascus13, Waki08, Waki12, Waki13} for more details. Hence, it is necessary to apply FR without solving problems that consist of LMI.  To propose such an approach, we provide  necessary and sufficient conditions in Theorem \ref{nonsf} by using FR that (\ref{LMI}) and its dual (\ref{dualLMI}) are not strongly feasible. In the next section, we will apply Theorem \ref{nonsf} into  LMI problems obtained from $H_{\infty}$ state feedback control problems and its dual.  
\theo\label{nonsf}
Problem (\ref{LMI}) is not strongly feasible if and only if there exists a nonzero $\hat{Y}\in\mathbb{S}^n$ such that 
\begin{align}\label{primalNSF}
F_j\bullet \hat{Y} = 0 \ (j\in\mathcal{M}),  F_0\bullet \hat{Y} \ge 0 \mbox{ and } \hat{Y}\in\mathbb{S}^n_+.
\end{align}
 In particular, if $\hat{Y}$ satisfies $F_0\bullet \hat{Y} > 0$, then (\ref{LMI}) is infeasible. If $F_0\bullet \hat{Y} =0$, (\ref{LMI}) is equivalent to the following problem: 
 \begin{align}\label{LMI_fr}
\inf_{x\in\mathbb{R}^m, X\in\mathbb{S}^n} \left\{c^Tx : X = \sum_{j\in\mathcal{M}}x_j F_j - F_0, X\in\mathbb{S}^n_+\cap\{\hat{Y}\}^{\perp}\right\}, 
 \end{align}
 where $\{\hat{Y}\}^{\perp}$ denotes the subspace  $\{X\in\mathbb{S}^n : X\bullet\hat{Y} = 0\}$ of $\mathbb{S}^n$. 
  Similarly, (\ref{dualLMI}) is not strongly feasible if and only if there exists a nonzero $(\hat{x}, \hat{X})\in\mathbb{R}^m\times\mathbb{S}^n$ such that 
 \begin{align}\label{dualNSF}
 \hat{X}=\displaystyle\sum_{j\in\mathcal{M}}\hat{x}_jF_j, \hat{X}\in\mathbb{S}^n_+ \mbox{ and }c^T\hat{x} \le 0. 
\end{align} 
 If $\hat{x}$ satisfies $c^T\hat{x} < 0$, then (\ref{dualLMI}) is infeasible. If $c^T\hat{x} = 0$, (\ref{dualLMI}) is equivalent to the following problem: 
  \begin{align}\label{dualLMI_fr}
  \sup_{Y\in\mathbb{S}^n} \left\{ F_0\bullet Y :  F_j\bullet Y = c_j \ (j\in\mathcal{M}), Y\in\mathbb{S}^n_+\cap\{\hat{X}\}^{\perp}
\right\}. 
 \end{align}
\etheo

The proof of  Theorem \ref{nonsf} are provided in \cite{Ramana97SIOPT, Pataki13, Waki13}. We give proofs of the if-part and the infeasibility in Appendix \ref{subsec:nonsf}.


\section{Conditions to be not strongly feasible}\label{sec:main}
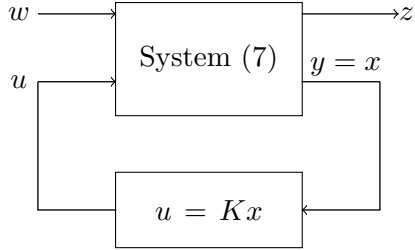
\begin{figure}[htbp]
\centering
\begin{tikzpicture}
    \tikzset{block/.style={draw,rectangle,text width=2.2cm, text centered, minimum height=1.5cm}};
    \node[block] (plant) {System (\ref{system})};
    \node  (w) at (-2.5, 0.6) {$w$};
   \node  (u) at (-2.5, -0.3) {$u$};
    \node  (z) at ( 2.6, 0.6) {$z$};
    \node  (y) at ( 1.8, -0.1) {$y=x$};
    \node[draw,rectangle,text width=2.2cm, text centered, minimum height=1.0cm] (control) at (0.0, -2.0) {$u=Kx$};
    \draw[->,line width=0.5pt] (-2.25, 0.6) -- (-1.22, 0.6);
    \draw[->,line width=0.5pt] (1.22, 0.6) -- (2.5, 0.6); 
    \draw[-,line width=0.5pt]  (1.22, -0.3) -- (2.25, -0.3);
    \draw[->,line width=0.5pt]  (2.25, -0.3) |- (1.22, -2.0);
    \draw[-,line width=0.5pt]  (-1.22, -2.0) -| (-2.25, -0.3);
    \draw[->,line width=0.5pt] (-2.25, -0.3) -- (-1.22, -0.3);
\end{tikzpicture}
\caption{The block diagram of  the closed loop system (\ref{clsystem})}
\label{fig:sfb}
\end{figure}

In this section, we give an interpretation of Theorem \ref{nonsf} in terms of control systems. 
We deal with $H_\infty$ state feedback control problem of a generalized plant. The generalized plant is given by 
\begin{equation}\label{system}
\left\{
\begin{array}{lclll}
\dot{x} &=& Ax + B_1 w + B_2 u\\
z &=& C_1x + D_{11}w +D_{12}u,
\end{array}
\right. 
\end{equation}
where $x\in\mathbb{R}^n$, $w\in\mathbb{R}^{m_1}$, $u\in\mathbb{R}^{m_2}$, $z\in\mathbb{R}^{p_1}$ 
and the matrices have the compatible dimensions. Let $u(t) = Kx(t)$ for a state feedback gain $K\in\mathbb{R}^{m_2\times n}$ to apply a state feedback in (\ref{system}). Then the state space representation of the closed loop system $G_{cl}(s)$ is given by 
\begin{equation}
\label{clsystem}
\left\{
\begin{array}{lcl}
\dot{x} &=&  (A+B_2K)x + B_1w \\
z &=&  (C_1+D_{12}K)x + D_{11}w.
\end{array}
\right.
\end{equation}
Figure \ref{fig:sfb} displays the block diagram of the closed loop system (\ref{clsystem}).

The following fact  is well-known on  (\ref{clsystem}). 
\theo\label{fundamental}(See {\itshape e.g.}, \cite{Iwasaki94})
For closed loop system (\ref{clsystem}) and a given $\gamma >0$, the following are equivalent: 
\renewcommand{\labelenumi}{(\Alph{enumi})}
\begin{enumerate}
\item There exists a $K\in\mathbb{R}^{m_2\times n}$ such that $\|G_{cl}(s)\|_{\infty} < \gamma$ and $A+B_2K$ is Hurwitz stable. 
\item There exist $X\in\mathbb{S}^n_{++}$ and $K\in\mathbb{R}^{m_2\times n}$ such that
\begin{align}\label{positivity}
&-\begin{pmatrix}
\He((A+B_2K)X) &* & *\\
(C_1+D_{12}K)X & -\gamma I_{p_1} & *\\
B_1^T & D_{11}^T & -\gamma I_{m_1}
\end{pmatrix}\in\mathbb{S}^{N_0}_{++}. 
\end{align}
\end{enumerate}
\etheo
Here $N_0=n+p_1+m_1$ and $*$ in (\ref{positivity})  of Theorem \ref{fundamental} 
stands for the transpose of the lower triangular block part. 

To obtain a state feedback gain $K$ for minimizing $\|G_{cl}(s)\|_{\infty}$,  one  can use the following LMI formulation by applying the change of variables method with $Y = KX$ in (\ref{positivity}):  
\begin{align}
\label{LMI2}
&\left\{
\begin{array}{cl}
\inf_{\gamma, X, Y} & \gamma\\
\mbox{sub. to} & -\begin{pmatrix}
\He(AX+B_2Y) &* & *\\
C_1X+D_{12}Y &-\gamma I_{p_1} & *\\
B_1^T & D_{11}^T & -\gamma I_{m_1}
\end{pmatrix}\in\mathbb{S}^{N_0}_+, \\
& \gamma\in\mathbb{R}, X\in\mathbb{S}^n_+, Y\in\mathbb{R}^{m_2\times n}, \\
\end{array}
\right.
\end{align}
where $N_0=n+p_1+m_1$ and $*$ stands for the transpose of the lower triangular block part. 

Its dual is formulated as follows. We describe a way to obtain (\ref{LMI2}) from (\ref{dualLMI2}) in Appendix \ref{sec:dual}. 
\begin{align}
\label{dualLMI2}
&\left\{
\begin{array}{cl}
\sup
& 2(B_1^T\bullet Z_{31} + D_{11}^T\bullet Z_{32})\\
\mbox{sub. to} & I_{p_1}\bullet Z_{22} + I_{m_1}\bullet Z_{33} = 1, B_2^TZ_{11} + D_{12}^TZ_{21} = O_{m_2\times n}, \\
& \He(A^TZ_{11} + C_1^TZ_{21})\in\mathbb{S}^n_+,  \begin{pmatrix}
Z_{11}&Z_{21}^T & Z_{31}^T\\
Z_{21}&Z_{22} &Z_{32}^T\\
Z_{31}&Z_{32}&Z_{33}
\end{pmatrix}\in\mathbb{S}^{N_0}_+. 
\end{array}
\right.
\end{align}
Here $Z_{ij} \ (1\le j\le i\le 3)$ is the decision variable in (\ref{dualLMI2}). It should be noted that LMI problem (\ref{LMI2}) is obtained by the change of variables method. One can also obtain another LMI formulation by the  elimination of variables method. We also remark that we can show the same results as those in this section. 


\subsection{Condition to be not strongly feasible for LMI problem (\ref{LMI2})}
Applying (\ref{primalNSF}) in Theorem \ref{nonsf} to LMI problem (\ref{LMI2}), we can obtain the condition to be not strongly feasible for (\ref{LMI2}). Namely, Problem  (\ref{LMI2}) is not strongly feasible if and only if the following problem on $Z_{ij}$ has a solution:  
\begin{align}\label{PNSFsystem}
\left\{
\begin{array}{l}
2(B_1^T\bullet Z_{31} + D_{11}^T\bullet Z_{32})\ge 0, I_{p_1}\bullet Z_{22} + I_{m_1}\bullet Z_{33} = 0, \\
\He(A^TZ_{11} + C_1^TZ_{21})\in\mathbb{S}^n_+, 
B_2^TZ_{11} + D_{12}^TZ_{21} = O_{m_2\times n}, \\ \begin{pmatrix}
Z_{11}&Z_{21}^T & Z_{31}^T\\
Z_{21}&Z_{22} &Z_{32}^T\\
Z_{31}&Z_{32}&Z_{33}
\end{pmatrix}\in\mathbb{S}^{N_0}_+\setminus\{O_{N_0\times N_0}\}. 
\end{array}
\right.
\end{align}

It follows from the equality constraint $I_{p_1}\bullet Z_{22} + I_{m_1}\bullet Z_{33} = 0$ and the positive semidefiniteness of $Z_{22}$ and $Z_{33}$ that any solution in (\ref{PNSFsystem}) satisfies $Z_{22}=O_{p_1\times p_1}$ and $Z_{33} = O_{m_1\times m_1}$, and thus (\ref{PNSFsystem}) is equivalent to the following problem on $Z_{11}$: 
\begin{align}\label{PNSFsystem2}
\He(A^TZ_{11})\in\mathbb{S}^n_+, B_2^TZ_{11} = O_{m_2\times n} \mbox{ and }
Z_{11}\in\mathbb{S}^n_+. 
\end{align}
We obtain the following proposition from  (\ref{PNSFsystem2}) and provide the  proof in Appendix \ref{sec:propPNSF}. 
\prop\label{propPNSF}
Problem (\ref{PNSFsystem}) has a solution if and only if 
there exists $\lambda\in\overline{\mathbb{C}_+}$ such that 
\begin{equation}\label{stabrank}
\rank\begin{pmatrix}
A - \lambda I_n, B_2
\end{pmatrix} < n, 
\end{equation}
{\itshape i.e.} $(A, B_2)$ is not stabilizable. 
\eprop
We obtain the following condition to be not strongly feasible for (\ref{LMI2}) by combining Proposition \ref{propPNSF} with Theorem \ref{nonsf}. 
\theo\label{stabilizable}
LMI problem (\ref{LMI2}) is strongly feasible if and only if control system (\ref{system}) is stabilizable. 
\etheo

Although Theorem \ref{stabilizable} has been already known in {\itshape e.g.}, \cite{Iwasaki94}, this theorem shows that one also can  prove the fact by facial reduction. 

\subsection{Condition to be not strongly feasible for Problem (\ref{dualLMI2})}\label{sec:dualfeas}
Applying  (\ref{dualNSF}) in Theorem \ref{nonsf} to Problem (\ref{dualLMI2}), we can obtain the condition to be not strongly feasible for (\ref{dualLMI2}). Namely, (\ref{dualLMI2}) is not strongly feasible if and only if the following problem on $\gamma$, $X$ and $Y$ has a nonzero solution: 
\begin{align}\label{DNSFsystem}
&\left\{
\begin{array}{l}
 -\begin{pmatrix}
\He(AX + B_2Y) &* & *\\
C_1X + D_{12}Y & -\gamma I_{p_1} &*\\
O_{m_1\times n} & O_{m_1\times p_1} & -\gamma I_{m_1}
\end{pmatrix}\in\mathbb{S}^{N_0}_+, \\
 \gamma\le 0, X\in\mathbb{S}^n_+, Y\in\mathbb{R}^{m_2\times n}. 
\end{array}
\right. 
\end{align}
Clearly, there does not exist any nonzero solutions $X$ and $Y$ such that (\ref{DNSFsystem}) holds and $\gamma <0$. Substituting $\gamma = 0$ into (\ref{DNSFsystem}), we obtain the following problem on $\gamma$, $X$ and $Y$ which is equivalent to (\ref{DNSFsystem}): 
\begin{align}\label{DNSFsystem2}
\left\{
\begin{array}{l}
C_1X + D_{12}Y =O_{p_1 \times n}, -\He(AX+B_2Y)\in\mathbb{S}^n_+, \\
 \gamma= 0, X\in\mathbb{S}^n_+, Y\in\mathbb{R}^{m_2\times n}. 
\end{array}
\right. 
\end{align}

We obtain the following proposition from (\ref{DNSFsystem2}) and provide the proof in Appendix \ref{sec:propDNSF}. 
\prop\label{propDNSF}
Suppose that the matrix $(B_2^T, D_{12}^T)^T$
is full column rank and that there exists $\lambda\in\overline{\mathbb{C}_-}$ such that 
\begin{align}\label{rankinvz}
&\rank \begin{pmatrix}
A - \lambda I_n & B_2 \\
C_1 & D_{12}
\end{pmatrix} < n+m_2. 
\end{align}
Then problem (\ref{DNSFsystem2}) has a solution. Furthermore, for a solution $(\gamma, X, Y)$ of (\ref{DNSFsystem2}), if there exists a full column rank matrix $H\in\mathbb{R}^{n\times r}$ and $R\in\mathbb{R}^{p_1\times r}$ such that $X = HH^T \mbox{ and } Y = RH^T$, then there exists $\lambda\in\overline{\mathbb{C}_-}$ such that (\ref{rankinvz}) holds. 
\eprop
We remark that $\lambda\in\overline{\mathbb{C}_-}$ that satisfies (\ref{rankinvz}) is called the {\itshape (stable) invariant zero of (\ref{system})}. We obtain the following theorem from the above discussion. 

\theo\label{invariantzeros}
If  $D_{12}$ is not full column rank, dual problem (\ref{dualLMI2}) is not strongly feasible. 
When $D_{12}$ is column full rank, dual problem (\ref{dualLMI2}) is not strongly feasible if and only if the subsystem $C_1(sI-A)^{-1}B_2+D_{12}$ of  control system (\ref{system}) has invariant zeros in $\overline{\mathbb{C}_-}$
\etheo
\begin{proof}
When $D_{12}$ is not full column rank, we can construct a nonzero $(\gamma, X, Y)$ that satisfies (\ref{DNSFsystem2}). In fact, for simplicity, we assume $D_{12} = (O_{p_1\times r}, \hat{D}_{12})$, where $\hat{D}_{12}\in\mathbb{R}^{p_1\times (m_2-r)}$ is full column rank. Denoting $B_2 = (B_{21}, B_{22})$, where $B_{21}\in\mathbb{R}^{n\times r}$, $B_{22}\in\mathbb{R}^{n\times (m_2-r)}$, we see that the following $(\gamma, X, Y)$ satisfies (\ref{DNSFsystem2}):
\[
\gamma = 0, X = O_{n\times n}, Y=\begin{pmatrix}
-B_{21}^T/2\\
O_{(m_2-r)\times n}
\end{pmatrix}. 
\]
It follows from Theorem \ref{nonsf} that dual problem (\ref{dualLMI2}) is not strongly feasible. 

If $D_{12}$ is full column rank, then there exists a full column rank matrix $H\in\mathbb{R}^{n\times r}$ such that $X=HH^T$ and $Y=RH^T$ for some $R\in\mathbb{R}^{p_1\times r}$.  In fact, we obtain $Y = -(D_{12}^TD_{12})^{-1}D_{12}^TC_1 X$ and $R = -(D_{12}^TD_{12})^{-1}D_{12}^TC_1H$ 
from $C_1X + D_{12}Y = O_{p_1\times n}$. Hence the desired result follows from Proposition \ref{propDNSF}. 
\end{proof}

\section{Reduction of LMI problem (\ref{LMI2}) and closed loop system (\ref{clsystem}) by using invariant zeros in $\overline{\mathbb{C}_-}$}\label{sec:reduction}

We present a reduction of LMI problem (\ref{LMI2}) to get rid of the numerical difficulty in solving it. The stabilizability in (\ref{system}) is a natural assumption for designing the state feedback gain $K$, whereas the existence of invariant zeros in $\overline{\mathbb{C}_-}$ in (\ref{system}) are not taken care with.  Dual problem (\ref{dualLMI2}) is not strongly feasible under the existence of invariant zeros in $\overline{\mathbb{C}_-}$ in (\ref{system}), and thus (\ref{dualLMI_fr}) in Theorem \ref{nonsf} is available to (\ref{dualLMI2}). The reduction that we present in this subsection is based of facial reduction and consists of the null vectors associated with invariant zeros in $\overline{\mathbb{C}_-}$. 

We show in this section that the resulting LMI problem is obtained from a subsystem of the closed loop system (\ref{clsystem}) and that the subsystem can be also obtained by applying the transformation with the non-singular matrix into (\ref{clsystem}).  In addition, we give a way to retrieve a state feedback gain $K$ which attains the optimal performance index from a state feedback gain of the subsystem and the non-singular matrix. 

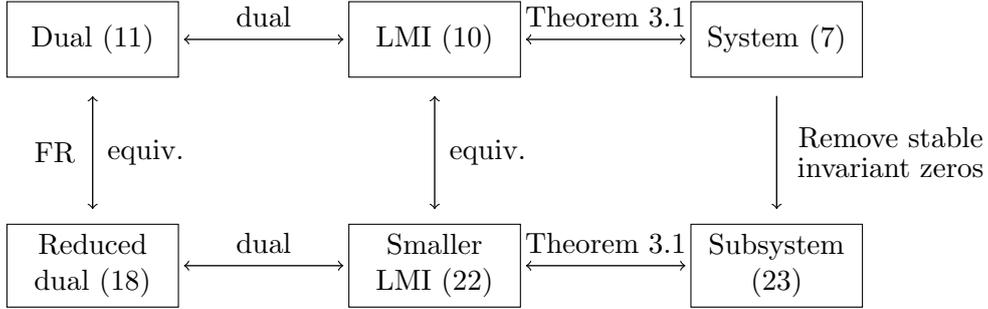
\begin{figure}[htpb]
\centering
\begin{tikzpicture}
    \tikzset{block/.style={draw,rectangle,text width=2.0cm, text centered, minimum height=1cm}};
    \node[block] (Dual) at (-4.5, 0.0) {Dual (\ref{dualLMI2})};
    \node[block] (LMI) {LMI (\ref{LMI2})};
    \node[block] (System) at (4.5, 0.0) {System (\ref{system})};
    \node[block] (RDual) at (-4.5, -3.0) {Reduced dual (\ref{dualLMI2_fr})};
    \node[block] (RLMI) at (0.0, -3.0) {Smaller LMI (\ref{LMI3})};
    \node[block] (SSystem) at (4.5, -3.0) {Subsystem (\ref{system2})};
    \node (msg1) at (-2.25, 0.3) {dual};
    \node (msg2) at (2.25, 0.3) {Theorem \ref{fundamental}};
    \node (msg8) at (2.25, -2.7) {Theorem \ref{fundamental}};
    \node (msg3) at (-2.25, -2.7) {dual};
    \node (msg4) at (1.5, -3.0) {};
    \node (msg5) at (-5.0, -1.5) {FR};
    \node (msg6) at (-3.8, -1.5) {equiv.};
    \node (msg7) at (0.7, -1.5) {equiv.};
    \node (msg9) at (6.0, -1.3) {Remove stable};
    \node (msg10) at (6.0, -1.7) {invariant zeros};
    \draw[<->,line width=0.5pt] (-3.3, 0.0) -- (-1.2, 0.0);
    \draw[<->,line width=0.5pt] (1.2, 0.0) -- (3.3, 0.0);
        \draw[<->,line width=0.5pt] (-3.3, -3.0) -- (-1.2, -3.0);
    \draw[<->,line width=0.5pt] (1.2, -3.0) -- (3.3, -3.0);
    \draw[<->,line width=0.5pt] (-4.5, -0.75) -- (-4.5, -2.25);
    \draw[<->,line width=0.5pt] (0.0, -0.75) -- (0.0, -2.25);
    \draw[->,line width=0.5pt] (4.5, -0.75) -- (4.5, -2.25);
\end{tikzpicture}
\caption{Relationship among control systems, LMIs and their duals in subsections \ref{subsec:reductionLMI} and \ref{subsec:reductionCLS}}
\label{sfb}
\end{figure}

In this section, we present reductions on LMI (\ref{LMI2}) and control system (\ref{system}) in subsections \ref{subsec:reductionLMI} and \ref{subsec:reductionCLS}. Figure \ref{sfb} displays the relationship among control systems, LMIs and their duals in subsections \ref{subsec:reductionLMI} and \ref{subsec:reductionCLS}. In subsection \ref{subsec:d12}, we deal with the case when $D_{12}$ in (\ref{system}) is not full column rank. 

Let $\lambda_j\in\overline{\mathbb{C}_-}$ be a invariant zeros in $\overline{\mathbb{C}_-}$ in (\ref{system}) for all $j=1, \ldots, r$, and denote the associated null vector by $(\eta_j, \xi_j)\in\mathbb{C}^{n}\times\mathbb{C}^{m_2}$, {\itshape i.e.}, for all $j=1, \ldots, r$, we have 
\[
\begin{pmatrix}
A -\lambda_j I_n & B_2\\
C_1 & D_{12}
\end{pmatrix}\begin{pmatrix}
\eta_j\\
\xi_j
\end{pmatrix} = 0. 
\]
We impose the following assumptions for simplicity: 
\assm
\begin{enumerate}[label=(A\arabic*)]
\item\label{A1} $\eta_1, \ldots, \eta_r$ are linearly independent, 
\item\label{A2} all stable invariant zeros $\lambda_j$ are real, and
\item\label{A3} no invariant zeros on the imaginary axis. 
\end{enumerate}
\eassm
In particular, Assumption \ref{A2} implies that both $\eta_j$ and $\xi_j$ are real for all $j=1, \ldots, r$. We remark that even if either Assumptions \ref{A1} or \ref{A2} fails, the discussion in Section \ref{sec:reduction} holds by applying a little technical manner. In contrast, it is more difficult to analyze system (\ref{system}) and the resulting LMI problem in which Assumption \ref{A3} fails. The analysis in the case is involved in future work. 

\subsection{Reduction of LMI problem (\ref{LMI2})}\label{subsec:reductionLMI}
It follows from the proof of Proposition \ref{propDNSF} that the following $(\hat{\gamma}, \hat{X}, \hat{Y})$ is nonzero solution of (\ref{DNSFsystem2}):
\[
\hat{\gamma} = 0, \hat{X} = \sum_{j=1}^r\eta_j\eta_j^T , \hat{Y} = \sum_{j=1}^r\xi_j\eta_j^T. 
\]
We define $\hat{H} = (\eta_1, \ldots, \eta_r)\in\mathbb{R}^{n\times r}$ and $\hat{R}=(\xi_1, \ldots, \xi_r)\in\mathbb{R}^{m_2\times r}$. $\hat{H}$ is full column rank and, $\hat{X}=\hat{H}\hat{H}^T$ and $\hat{Y}=\hat{R}\hat{H}^T$.  Let $\hat{W}:=-\He(A\hat{X}+B_2\hat{Y})$. It follows from Lemma \ref{techlemma} that we have $A\hat{H} + B_2\hat{R} =\hat{H}\mbox{Diag}(\lambda_1, \ldots, \lambda_r)$, and thus we have 
$\hat{W} = -2\hat{H}\mbox{Diag}(\lambda_1, \ldots, \lambda_r)\hat{H}^T$. 
Here $\mbox{Diag}(x_1, \ldots, x_r)$ stands for the diagonal matrix with $x_1, \ldots, x_r$. In addition, it follows from Assumption \ref{A3} that the rank of $\hat{W}$ is $r$. 

 By applying (\ref{dualLMI_fr}) in Theorem \ref{nonsf} to dual problem (\ref{dualLMI2}), we obtain the  following problem: 
\begin{equation}
\label{dualLMI2_fr}
\left\{
\begin{array}{cl}
\sup
& 2(B_1^T\bullet Z_{31} + D_{11}^T\bullet Z_{32})\\
\mbox{sub. to} & I_{p_1}\bullet Z_{22} + I_{m_1}\bullet Z_{33} = 1, B_2^TZ_{11} + D_{12}^TZ_{21} = O_{m_2\times n}, \\
& \He(A^TZ_{11} + C_1^TZ_{21}) \in\mathbb{S}^n_+\cap\{\hat{X}\}^{\perp}, \\
& \begin{pmatrix}
Z_{11}&Z_{21}^T & Z_{31}^T\\
Z_{21}&Z_{22} &Z_{32}^T\\
Z_{31}&Z_{32}&Z_{33}
\end{pmatrix}\in\mathbb{S}^{N_0}_+\cap\left\{
\begin{pmatrix}
\hat{W} &O&O\\
O&O&O\\
O&O&O
\end{pmatrix}
\right\}^{\perp}. 
\end{array}
\right.
\end{equation}
Problem (\ref{dualLMI2_fr}) is equivalent to dual (\ref{dualLMI2}) of  the original LMI problem (\ref{LMI2}).  From the constraints in (\ref{dualLMI2_fr}), we have 
$\hat{H}^T\He(A^TZ_{11} + C_1^TZ_{21}) \hat{H} = O_{r\times r}$ and $Z_{k1}\hat{H} = O$ for $k=1, 2, 3$
 for all solutions of (\ref{dualLMI2_fr}).  
 
One can construct $J\in\mathbb{R}^{n\times (n-r)}$ so that  the matrix $T = (\hat{H}, J)$ is non-singular since $\hat{H}$ is full column rank. From these constraints in (\ref{dualLMI2_fr}), for any solution in (\ref{dualLMI2_fr}), we have
\begin{align}
\nonumber
& \begin{pmatrix}
T^T & &\\
& I_{p_1} & \\
&&I_{m_1} 
\end{pmatrix}
\left(\begin{array}{c|cc}
Z_{11}&Z_{21}^T & Z_{31}^T\\
\hline
Z_{21}&Z_{22} &Z_{32}^T\\
Z_{31}&Z_{32}&Z_{33}
\end{array} 
\right)
\begin{pmatrix}
T& &\\
& I_{p_1} & \\
&&I_{m_1} 
\end{pmatrix}\\
\label{reductionL}
&=\left(\begin{array}{cc|cc}
O_{r\times r} & O_{r\times (n-r)} & O_{r\times p_1} & O_{r\times m_1} \\
O_{(n-r)\times r}&\tilde{Z}_{11} & \tilde{Z}_{21}^T & \tilde{Z}_{31}^T\\
\hline
O_{p_1\times r} & \tilde{Z}_{21} & Z_{22} & Z_{32}^T\\
O_{m_1\times r}&\tilde{Z}_{31} & Z_{32} & Z_{33}
\end{array} 
\right). 
\end{align}

We define coefficient matrices $\tilde{A}, \tilde{B}_i \ (i=1, 2)$ and $\tilde{C}_1$ by  
\begin{align}\label{tildeCoef}
T^{-1}AT &= \begin{pmatrix}
\tilde{A}_{11} & \tilde{A}_{12}\\
\tilde{A}_{21} & \tilde{A}_{22}
\end{pmatrix}, T^{-1}B_i = \begin{pmatrix}
\tilde{B}_{i1}\\
\tilde{B}_{i2}
\end{pmatrix}, 
C_1T = \begin{pmatrix}\tilde{C}_{11}, \tilde{C}_{12}\end{pmatrix}. 
\end{align}

Then, (\ref{dualLMI2_fr}) can be reformulated as follows: 
\begin{equation}
\label{dualLMI2_fr2}
\left\{
\begin{array}{cl}
\sup
& 2(\tilde{B}_{12}^T\bullet \tilde{Z}_{31} + D_{11}^T\bullet Z_{32})\\
\mbox{sub. to} & I_{p_1}\bullet Z_{22} + I_{m_1}\bullet Z_{33} = 1, \begin{pmatrix}
\tilde{A}_{21}^T\\
\tilde{B}_{22}^T
\end{pmatrix}\tilde{Z}_{11} + 
\begin{pmatrix}
\tilde{C}_{11}^T\\
D_{12}^T
\end{pmatrix}\tilde{Z}_{21} = O, \\
& \He(\tilde{A}_{22}^T\tilde{Z}_{11} + \tilde{C}_{12}^T\tilde{Z}_{21}) \in\mathbb{S}^{n-r}_+, \begin{pmatrix}
\tilde{Z}_{11}&\tilde{Z}_{21}^T & \tilde{Z}_{31}^T\\
\tilde{Z}_{21}&Z_{22} &Z_{32}^T\\
\tilde{Z}_{31}&Z_{32}&Z_{33}
\end{pmatrix}\in\mathbb{S}^{N}_+, 
\end{array}
\right.
\end{equation}
where $N=n-r+p_1+m_1$.  In Appendix \ref{appendix:computation}, we give a detail how to obtain (\ref{dualLMI2_fr2}). 

Denote $T^{-1}= 
\begin{pmatrix}
U_1\\
U_2
\end{pmatrix}$. 
Since we have $A\hat{H} + B_2\hat{R} = \hat{H}\mbox{Diag}(\lambda_1, \ldots, \lambda_r)$, $C_1\hat{H} + D_{12}\hat{R} = O_{p_1\times r}$ and  $U_2\hat{H} = O_{(n-r)\times r}$, we obtain 
$\tilde{A}_{21} = \tilde{B}_{22}\hat{R}$ and $\tilde{C}_{11} = -D_{12}\hat{R}$, and thus the equality constraint $\tilde{A}_{21}^T\tilde{Z}_{11}+\tilde{C}_{11}^T\tilde{Z}_{21}=O_{r\times (n-r)}$ in (\ref{dualLMI2_fr2}) is redundant. 

The primal LMI problem of (\ref{dualLMI2_fr2}) is 
\begin{equation}
\label{LMI3}
\left\{
\begin{array}{cl}
\inf_{\tilde{\gamma}, \tilde{X}, \tilde{Y}} & \tilde{\gamma}\\
\mbox{sub. to} & -\begin{pmatrix}
\He(\tilde{A}_{22}\tilde{X} + \tilde{B}_{22}\tilde{Y})& * & *\\
\tilde{C}_{12}\tilde{X} + D_{12}\tilde{Y} & -\tilde{\gamma} I_{p_1}&*\\
\tilde{B}_{12}^T & D_{11}^T & -\tilde{\gamma} I_{m_1}
\end{pmatrix}\in\mathbb{S}^{N}_+, \\
& \tilde{\gamma}\in\mathbb{R}, \tilde{X}\in\mathbb{S}^{n-r}_+, \tilde{Y}\in\mathbb{R}^{m_2\times (n-r)}.  
\end{array}
\right.
\end{equation}
We remark that the optimal value $\tilde{\gamma}^*$ of (\ref{LMI3}) is equivalent to the optimal value $\gamma^*$ of the original LMI problem (\ref{LMI2}). 

\subsection{Reduction of closed loop system (\ref{clsystem})}\label{subsec:reductionCLS}

We discuss the reduction of control system (\ref{system}) via the reduction discussed in the previous subsection. To this end we focus on 
LMI problem (\ref{LMI3}). This is obtained  by applying Theorem \ref{fundamental} into
from the system whose state space representation is 
\begin{equation}\label{system2}
\left\{
\begin{array}{ccl}
\dot{\tilde{x}}_2 &=& \tilde{A}_{22}\tilde{x}_2 +\tilde{B}_{12}w + \tilde{B}_{22}u\\
z &=& \tilde{C}_{12}\tilde{x}_2 + D_{11}w + D_{12}u, 
\end{array}
\right.
\end{equation}
where the dimension of $\tilde{x}_2$  is $n-r$. 
Let $\tilde{K}\in\mathbb{R}^{m_2\times (n-r)}$ be a state feedback gain for (\ref{system2}) whose performance index is the optimal value $\tilde{\gamma}^*$ of (\ref{LMI3}), {\itshape i.e.}, it is equal to the optimal value $\gamma^*$ of the original LMI problem (\ref{LMI2}).  Then, the closed loop system is provided by 
\begin{equation}\label{clsubsystem2}
\left\{
\begin{array}{ccl}
\dot{\tilde{x}}_2 &=&(\tilde{A}_{22}+\tilde{B}_{22}\tilde{K})\tilde{x}_2 +\tilde{B}_{12}w\\
z &=& (\tilde{C}_{12}+D_{12}\tilde{K})\tilde{x}_2 + D_{11}w. 
\end{array}
\right.
\end{equation}

We show that the  state feedback gain $K$ for (\ref{system}), which is defined by $K = (\hat{R}, \tilde{K})T^{-1}$, attains the performance index $\tilde{\gamma}^*$. To this end, we apply the transformation with $T$ into (\ref{clsystem}). Denote $\tilde{x}(t) = T^{-1}x(t)$, then we obtain 
\begin{equation}
\label{clsystem3}
\left\{
\begin{array}{ccl}
\dot{\tilde{x}} &=& T^{-1}(A+B_2K)T\tilde{x} +T^{-1}B_1w\\
z &=& (C_1 + D_{12}K)T\tilde{x} + D_{11}w.
\end{array}
\right. 
\end{equation}
Recall $U_1T = (I_r, O_{r\times (n-r)})$ and $U_2\hat{H}=O_{(n-r)\times r}$. It follows that we have 
\begin{align*}
T^{-1}(A+B_2K)T
&= \begin{pmatrix}
\mbox{Diag}(\lambda_1, \ldots, \lambda_r)& U_1\left(AJ+B_2\tilde{K}\right)\\
O_{(n-r)\times r} & U_2\left(AJ+B_2\tilde{K}\right)
\end{pmatrix}, \\
(C_1+D_{12}K)T &= \begin{pmatrix}
O_{p_1\times r} & C_1J+D_{12}\tilde{K}
\end{pmatrix}.  
\end{align*}

Let $\tilde{x}(t)=(\tilde{x}_1(t)^T, \tilde{x}_2(t)^T)^T$. All $\lambda_j$  are 
unobservable modes, {\itshape i.e.}, $\tilde{x}_1$ is unobservable. Since $\tilde{x}_2$ is controllable and observable, performance index of the following subsystem (\ref{clsubsystem}) of (\ref{clsystem3}) is equivalent to the performance index of (\ref{clsystem3}): 
\begin{equation}
\label{clsubsystem}
\left\{
\begin{array}{ccl}
\dot{\tilde{x}}_2 &=&U_2(AJ+B_2\tilde{K})
\tilde{x}_2 + U_2B_1w\\
z &=& (C_1J+D_{12}\tilde{K})\tilde{x}_2 + D_{11}w. 
\end{array}
\right.
\end{equation}
Since we have $\tilde{A}_{22} = U_2AJ$, $\tilde{B}_{12}=U_2B_1$, $\tilde{B}_{22}=U_2B_2$ and $\tilde{C}_{12} = C_1J$, (\ref{clsubsystem}) is equivalent to (\ref{clsubsystem2}). Therefore, the performance index of the gain $K= (\hat{R}, \tilde{K})T^{-1}$ is the optimal value ${\gamma}^*$ of LMI problem (\ref{LMI2}). As shown above, the part of state feedback gain defined by $\hat{R}$  eliminates the invariant zeros by pole-zero cancellation. Theorem \ref{invariantzeros} ensures that this pole-zero cancellation does not affect the optimal value $\gamma$. 

\subsection{Reduction in the case where $D_{12}$ is not full column rank}\label{subsec:d12}
We consider the case where $D_{12}$ is not full column rank. For simplicity, we assume $D_{12} = (O_{p_1\times r}, \hat{D}_{12})$, where $\hat{D}_{12}\in\mathbb{R}^{p_1\times (m_2-r)}$ is full column rank. Denote $B_2=(B_{21}, B_{22})$, where $B_{21}\in\mathbb{R}^{n\times r}$ and $B_{22}\in\mathbb{R}^{n\times (m_2-r)}$. Let $(\gamma, X, Y)$ be as in the proof of Theorem \ref{invariantzeros}, {\itshape i.e.} 
$\gamma = 0$, $X = O_{n\times n}$ and 
\[
Y = \begin{pmatrix}
-B_{21}^T\\
O_{(m_2-r)\times n}
\end{pmatrix}. 
\]
Then $(\gamma, X, Y)$ satisfies (\ref{DNSFsystem2}). We define $\hat{W} = -\He(AX + B_2Y) = B_{21}B_{21}^T$, and thus the reduced dual problem can be formulated as follows: 
\begin{equation}
\label{dualLMI2_d12_fr}
\left\{
\begin{array}{cl}
\sup
& 2(B_1^T\bullet Z_{31} + D_{11}^T\bullet Z_{32})\\
\mbox{sub. to} & I_{p_1}\bullet Z_{22} + I_{m_1}\bullet Z_{33} = 1, B_2^TZ_{11} + D_{12}^TZ_{21} = O_{m_2\times n}, \\
& \He(A^TZ_{11} + C_1^TZ_{21}) \in\mathbb{S}^n_+, \\
& \begin{pmatrix}
Z_{11}&Z_{21}^T & Z_{31}^T\\
Z_{21}&Z_{22} &Z_{32}^T\\
Z_{31}&Z_{32}&Z_{33}
\end{pmatrix}\in\mathbb{S}^{N_0}_+\cap\left\{
\begin{pmatrix}
\hat{W} &O&O\\
O&O&O\\
O&O&O
\end{pmatrix}
\right\}^{\perp}. 
\end{array}
\right.
\end{equation}

For (\ref{dualLMI2_d12_fr}), it follows from the last constraint that we have $Z_{k1}B_{21} = O$ \ $(k=1, 2, 3)$. Hence the equality constraint $B_2^TZ_{11} + D_{12}^TZ_{21} = O_{m_2\times n}$ is equivalent to the equality constraint $B_{22}^TZ_{11} + \hat{D}_{12}^TZ_{21} = O_{(m_2-r)\times n}$. In addition, there exists $L\in\mathbb{R}^{n\times (n-r)}$ such that $T:=(B_{21}, L)$ is non-singular. Then  (\ref{reductionL}) holds for any feasible solution $Z_{ij}$. We define matrices $\tilde{A}_{ij}$, $\tilde{B}_{1j}$ and $\tilde{C}_{1j}$ as in (\ref{tildeCoef}). Then we have 
\begin{align*}
&T^T\He(A^TZ_{11} + C_1^TZ_{21})T \\
&= \He\left(
\begin{pmatrix}
\tilde{A}^T_{11} & \tilde{A}^T_{21}\\ 
\tilde{A}^T_{12} & \tilde{A}^T_{22}
\end{pmatrix}\begin{pmatrix}
O_{r\times r} & O_{r\times(n-r)}\\
O_{(n-r)\times r} & \tilde{Z}_{11}
\end{pmatrix} + \begin{pmatrix}
\tilde{C}_{11}^T\\
\tilde{C}_{12}^T
\end{pmatrix}\begin{pmatrix}
O_{p_1\times r} & \tilde{Z}_{21}
\end{pmatrix}
\right)\\
&=\He\left(
\begin{pmatrix}
O_{r\times r} & \tilde{A}^T_{21}\tilde{Z}_{11}+\tilde{C}_{11}^T\tilde{Z}_{21}\\
O_{(n-r)\times r} & \tilde{A}^T_{22}\tilde{Z}_{11}+\tilde{C}_{12}^T\tilde{Z}_{21}
\end{pmatrix}
\right). 
\end{align*}
Thus any feasible solution of (\ref{dualLMI2_d12_fr}) satisfies
\[
\tilde{A}^T_{21}\tilde{Z}_{11}+\tilde{C}_{11}^T\tilde{Z}_{21} = O_{r\times (n-r)}\mbox{ and } \tilde{A}^T_{22}\tilde{Z}_{11}+\tilde{C}_{12}^T\tilde{Z}_{21}\in\mathbb{S}^{n-r}_+. 
\]
Therefore we obtain the following dual problem which is equivalent to (\ref{dualLMI2_d12_fr}): 
\begin{align}
\label{Dual4}
&\left\{
\begin{array}{cl}
\sup& 2(\tilde{B}_{12}^T\bullet \tilde{Z}_{31} + D_{11}^T\bullet Z_{32})\\
\mbox{sub. to} & I_{p_1}\bullet Z_{22} + I_{m_1}\bullet Z_{33} = 1, E_1^T\tilde{Z}_{11} + 
E_2^T\tilde{Z}_{21} = O_{m_2\times (n-r)},\\
&\He\left(\tilde{A}_{22}^T\tilde{Z}_{11} + \tilde{C}_{12}^T\tilde{Z}_{21}\right) \in\mathbb{S}^{n-r}_+,  \begin{pmatrix}
\tilde{Z}_{11}&\tilde{Z}_{21}^T & \tilde{Z}_{31}^T\\
\tilde{Z}_{21}&Z_{22} &Z_{32}^T\\
\tilde{Z}_{31}&Z_{32}&Z_{33}
\end{pmatrix}\in\mathbb{S}^{N}_+, 
\end{array}
\right. 
\end{align}
where $N=n-r+p_1+m_1$, $E_1 := (\tilde{A}_{21}, \tilde{B}_{222})\in\mathbb{R}^{n\times m_2}$, $E_2 := (\tilde{C}_{11}, \hat{D}_{12})\in\mathbb{R}^{p_1\times m_2}$ and we decompose $B_{22}$ as follows: 
 \[
 T^{-1} B_{22} = \begin{pmatrix}
 \tilde{B}_{221}\\
 \tilde{B}_{222}
 \end{pmatrix}. 
 \]
 
 LMI problem which corresponds to (\ref{dualLMI2_d12_fr}) can be formulated as follows:
 \begin{align}
\label{LMI4}
&\left\{
\begin{array}{cl}
\inf & \gamma\\
\mbox{sub. to} & \gamma\in\mathbb{R}, X\in\mathbb{S}^{n-r}_+, Y\in\mathbb{R}^{m_2\times (n-r)}, \\
& -\begin{pmatrix}
\He(\tilde{A}_{22}X + E_1Y) &* & *\\
\tilde{C}_{12}X + E_2Y & -\gamma I_{p_1} & *\\
\tilde{B}_{12}^T & D_{11}^T & -\gamma I_{m_1}
\end{pmatrix}\in\mathbb{S}^{N}_+
\end{array}
\right. 
\end{align}
 
 A closed loop system which corresponds to LMI (\ref{LMI4}) is provided as follows: 
\begin{equation}
\label{rsss}
\left\{
\begin{array}{ccl}
\dot{\tilde{x}}(t) &=&  \left(\tilde{A}_{22}+\left(\tilde{A}_{21}, \tilde{B}_{222}\right)\tilde{K}\right)\tilde{x}(t) + \tilde{B}_{12}w(t)  \\
z(t) &=&  \left(\tilde{C}_{12}+\left(\tilde{C}_{11}, \hat{D}_{12}\right)\tilde{K}\right)\tilde{x}(t) + D_{11}w(t), 
\end{array}
\right.
\end{equation}
where $\tilde{x}\in \mathbb{R}^{n-r}$ is a state variable and $\tilde{K}\in\mathbb{R}^{m_2\times (n-r)}$ is a state feedback gain.

\begin{figure}[htbp]
\centering
\begin{tikzpicture}
    \tikzset{block/.style={draw,rectangle,text width=2.2cm, text centered, minimum height=1.5cm}};
    \node[block] (plant) {System (\ref{system})};
    \node  (w) at (-4.5, 0.6) {$w$};
    \node  (u1) at (-1.5, 0.2) {$u_1$};
    \node  (u2) at (-1.5, -0.4) {$u_2$};
    \node   (ut) at (-3.4, 0.25) {$\tilde{u}_1$};
    \node  (z) at ( 2.7, 0.6) {$z$};
    \node  (y) at ( 1.8, -0.3) {$y=x$};
    \node[circle, draw, fill=black, inner sep=0.75pt]  (c) at (-3.8, -0.3) {};
    \node[draw,rectangle,text width=1cm, text centered, minimum height=0.25cm] (diff) at (-2.5, 0.0) {$s+\alpha$};
    \node[draw,rectangle,text width=2.2cm, text centered, minimum height=1.0cm] (control) at (0.0, -2.0) {$\begin{pmatrix}\tilde{u}_1\\ u_2\end{pmatrix}=K_dx$};
    \draw[->,line width=0.5pt] (-4.25, 0.6) -- (-1.22, 0.6);
    \draw[->,line width=0.5pt]  (-3.8, -0.6) -- (-1.22, -0.6);
    \draw[->,line width=0.5pt] (1.22, 0.6) -- (2.5, 0.6); 
    \draw[-,line width=0.5pt]  (1.22, -0.6) -- (2.25, -0.6);
    \draw[->,line width=0.5pt]  (2.25, -0.6) |- (1.22, -2.0);
    \draw[-,line width=0.5pt]  (-1.22, -2.0) -| (-4.25, -0.3);
    \draw[->,line width=0.5pt] (-1.88, 0.0) -- (-1.22, 0.0);
    \draw[->,line width=0.5pt] (-3.8, -0.6) |- (-3.14, 0.0);
    \draw[-,line width=0.5pt] (-4.25, -0.3) -- (-3.8, -0.3);
\end{tikzpicture}
\caption{The block diagram of  the closed loop system (\ref{diffcl})}
\label{fig:sfb2}
\end{figure}
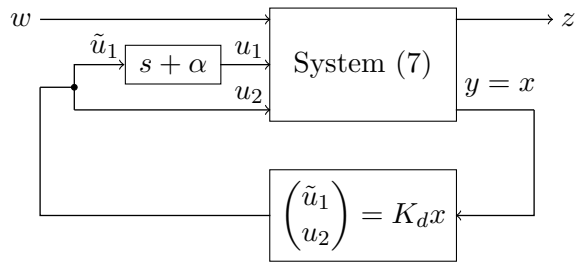

We provide an interpretation of  the performance index $\gamma$ of (\ref{rsss}). The original system (\ref{system}) can be described as follows:
\[
\left\{
\begin{array}{ccl}
\dot{x} &=& Ax + B_1 w + \begin{pmatrix}
B_{21} & B_{22}
\end{pmatrix}\begin{pmatrix}
u_1\\
u_2
\end{pmatrix}\\
z &=& C_1 x + D_{11}w + \begin{pmatrix}
O_{p_1\times r} & \hat{D}_{12}
\end{pmatrix}\begin{pmatrix}
u_1\\
u_2
\end{pmatrix}. 
\end{array}
\right.
\]
We consider the following control system with a differentiator added in the above system. 
\begin{equation}\label{diffsss}
\left\{
\begin{array}{ccl}
\dot{x} &=& Ax + B_1 w + (A+\alpha I_n) B_{21}\tilde{u}_1 + B_{22}u_2\\
z &=& C_1x  + D_{11}w +C_1B_{21} \tilde{u}_1 + \hat{D}_{12}u_2. 
\end{array}
\right.
\end{equation}
Here $\alpha >0$. We consider the state feedback controller $\tilde{u}_1 = K_{d1}x$ and $u_2 = K_{d2}x$, where $K_{d1}\in\mathbb{R}^{r\times n}$ and $K_{d2}\in\mathbb{R}^{(m_2-r)\times n}$. The closed loop system shown in Figure \ref{fig:sfb2} is formulated as follows: 
\begin{equation}\label{diffcl}
\left\{
\begin{array}{ccl}
\dot{x} &=& (A + (A+\alpha I_n)B_{21}K_{d1} + B_{22}K_{d2})x + B_1 w\\
z &=& (C_1 +C_1B_{21} K_{d1}+ \hat{D}_{12}K_{d2}) x + D_{11}w. 
\end{array}
\right.
\end{equation}

We remark that the closed loop system (\ref{diffcl}) has  invariant zeros in $\overline{\mathbb{C}_-}$. In fact, we have $AB_{21} +(A+\alpha I_n)B_{21}(-I_r) = -\alpha B_{21}$ and $C_1B_{21} + C_1B_{21}(-I_r) = O_{p_1\times r}$. Since the closed loop (\ref{diffcl})  has invariant zeros in $\overline{\mathbb{C}_-}$, its dual of  LMI problem obtained from (\ref{diffsss}) is not strongly feasible. Applying a similar discussion in subsection \ref{subsec:reductionLMI} and $T=(B_{21}, L)$ to the dual problem, we obtain the same LMI problem as (\ref{LMI4}). Hence the optimal performance index of (\ref{diffcl}) is equivalent to that of (\ref{clsystem}). 

We define the state feedback gain $K_d$ by 
\[
K_d = \begin{pmatrix}
K_{d1}\\
K_{d2}
\end{pmatrix} := \begin{pmatrix}
-I_r & \tilde{K}_{d1}\\
O_{(m_2-r)\times r} & \tilde{K}_{d2} 
\end{pmatrix}T^{-1}, \tilde{K}_d = \begin{pmatrix}
\tilde{K}_{d1}\\
\tilde{K}_{d2}
\end{pmatrix}
\]
where $\tilde{K}_{d1}\in\mathbb{R}^{r\times (n-r)}$ and $\tilde{K}_{d2}\in\mathbb{R}^{(m_2-r)\times (n-r)}$. 
Here $\tilde{K}_d$ is the optimal state feedback obtained by solving (\ref{LMI4}). 
By applying the transformation with $T$ into the closed system (\ref{diffcl}), we obtain
\begin{equation}\label{diffcl_trans}
\left\{
\begin{array}{ccl}
\dot{\tilde{x}} &=& T^{-1}(A + (A+\alpha I_n)B_{21}K_{d1} + B_{22}K_{d2})T\tilde{x}  + T^{-1}B_1 w\\
z &=& (C_1 +C_1B_{21} K_{d1}+ \hat{D}_{12}K_{d2}) T\tilde{x} + D_{11}w. 
\end{array}
\right.
\end{equation}
By direct computation, we have in (\ref{diffcl_trans}), 
\begin{align*}
T^{-1}(A+\alpha I_n)B_{21}K_{d1}T &=T^{-1}(A+\alpha I_n)B_{21}\begin{pmatrix}
-I_{r} & \tilde{K}_{d1}
\end{pmatrix} \\
&= \begin{pmatrix}
-\tilde{A}_{11} - \alpha I_r & (\tilde{A}_{11}+\alpha I_r)\tilde{K}_{d1}\\
-\tilde{A}_{21} & \tilde{A}_{21}\tilde{K}_{d1}
\end{pmatrix}, 
\\
T^{-1}B_{22}K_{d2}T &=\begin{pmatrix}
O_{n\times r} & T^{-1}B_{22}\tilde{K}_{d2}
\end{pmatrix}, \\
C_1B_{21} K_{d1}T&=
\begin{pmatrix}
-\tilde{C}_{11} & \tilde{C}_{11}\tilde{K}_{d1}
\end{pmatrix}, 
 \hat{D}_{12}K_{d2}T = \begin{pmatrix}
O_{p_1\times r} & \hat{D}_{12}\tilde{K}_{d2}
\end{pmatrix}. 
\end{align*}
From these direct computation, (\ref{diffcl_trans}) is equivalent to
\begin{equation}\label{diffcl_trans2}
\left\{
\begin{array}{ccl}
\begin{pmatrix}
\dot{\tilde{x}}_1\\
 \dot{\tilde{x}}_2
\end{pmatrix}
&=& \begin{pmatrix}
-\alpha I_r & \tilde{A}_{12} + (\tilde{A}_{11}+\alpha I_r)\tilde{K}_{d1} + \tilde{B}_{221}\tilde{K}_{d2}\\
O_{(n-r)\times r} & \tilde{A}_{22} + \tilde{A}_{21}\tilde{K}_{d1}+\tilde{B}_{222}\tilde{K}_{d2}
\end{pmatrix}
\begin{pmatrix}
\tilde{x}_1\\
\tilde{x}_2
\end{pmatrix}+ T^{-1}B_1 w\\
z &=& \begin{pmatrix} 
O_{p_1\times r} & \tilde{C}_{12} + \tilde{C}_{11}\tilde{K}_{d1} +  \hat{D}_{12}\tilde{K}_{d2}
\end{pmatrix}
\begin{pmatrix}
\tilde{x}_1\\
\tilde{x}_2
\end{pmatrix}+ D_{11}w. 
\end{array}
\right.
\end{equation}
We see form (\ref{diffcl_trans2}) that (\ref{rsss}) is the part of $\tilde{x}_2$. In addition, $\tilde{K}_d$ is the same as $\tilde{K}$ in (\ref{rsss}) since both the obtained LMI problems are the same. We remark that the part of $\tilde{x}_2$, $\tilde{K}_{d1}$,  $\tilde{K}_{d2}$ and $K_d$ are  also independent on $\alpha$ used in the differentiator. This fact implies that the same optimal performance index can be achieved
for systems (\ref{system}) and (\ref{diffsss}) by static feedback gains although the optimal feedback gains $K$ and $K_d$ are different. 

\section{Numerical experiment}\label{sec:experiment}
In this section, we compare numerical performance of reduced LMI problems (\ref{LMI3}) and (\ref{LMI4})  with the original LMI problem (\ref{LMI2}) obtained from the original system (\ref{system}).  We use 
SDPT3 to solve LMI problems by calling it from YALMIP \cite{Lofberg04}\footnote{We modified {\ttfamily computedimacs.m} in YALMIP to obtain more correct DIMACS error $\mbox{err}_3$ and $\mbox{err}_6$.} with default parameters. We set the stopping tolerance $\epsilon = 1.0\times 10^{-7}$.  

We use DIMACS errors to see the numerical performance, which is defined in \cite{Mittelmann03}. They measure the accuracy of solutions $(x, X)$ and $Y$ obtained by LMI software  for LMI problem (\ref{LMI}) and its dual (\ref{dualLMI}), and consists of six errors. We introduce the following {three} of DIMACS errors: 

\begin{align*}
\mbox{err}_1 &:= \frac{\displaystyle\sqrt{\sum_{j\in\mathcal{M}}\left(F_j\bullet Y-c_j\right)^2}}{1+\|c\|_1},  
\mbox{err}_5 := \frac{c^Tx - F_0\bullet Y}{1+|c^Tx|+|F_0\bullet Y|}, \mbox{err}_6 := \frac{X\bullet Y}{1+|c^Tx|+|F_0\bullet Y|}, 
\end{align*}
where $\|\cdot\|_1$ is the largest absolute value in the vector.  
In this manuscript, we call $\mbox{err}_1$, 
$\mbox{err}_5$ and  $\mbox{err}_6$ by {\itshape dual feasibility}, 
{\itshape relative gap} and {\itshape relative complementarity}, respectively. We omit information on $\mbox{err}_2$, $\mbox{err}_3$ and $\mbox{err}_4$  because they are always zeros for computed solutions.  
Since the stopping tolerance is $\epsilon = 1.0\times 10^{-7}$, the computed dual solution is not feasible when $\mbox{err}_1$ is bigger than $\epsilon$. In addition, when $\mbox{err}_5$ and/or $\mbox{err}_6$ is bigger than $\epsilon$, the computed solution is not optimal. 


We describe a way to generate control system (\ref{system}) that has given stable invariant zeros $\lambda$ in Algorithm \ref{algo}. In this numerical experiment, we remark that each element in coefficient matrices $A$, $C_1$, $B_1$, $B_2$, $D_{11}$ and $D_{12}$  is in $[-3, 3]$. 

\normalem
\begin{algorithm}[H] \label{algo}
\KwIn{$(n, p_1, m_1, m_2)\in\mathbb{N}^4$, $r\in\mathbb{N}$ and stable invariant zeros $\lambda_1, \ldots, \lambda_r\in\overline{\mathbb{C}_-}$}
\KwOut{$A\in\mathbb{R}^n$, $C_1\in\mathbb{R}^{p_1\times n}$, $B_1\in\mathbb{R}^{n\times m_1}$, $B_2\in\mathbb{R}^{n\times m_2}$, $D_{11}\in\mathbb{R}^{p_1\times m_1}$ and $D_{12}\in\mathbb{R}^{p_1\times m_2}$}
$A$, $C_1$, $B_1$, $B_2$, $D_{11}$ and $D_{12}$ are randomly generated\;
\tcc{We denote the $i$th column of $A$ and $C_1$ by $a_i$ and $c_i$, respectively. }
\tcc{Let $e_j\in\mathbb{R}^r$ be the $r$-dimensional $j$th unit vector. }
\For{$j\rightarrow 1$ \KwTo $r$}{
$\begin{pmatrix}
\hat{\eta}_j\\
\xi_j
\end{pmatrix}\in\mathbb{R}^{(n-r)+m_2}$ is randomly generated\; 
$
\begin{pmatrix}
v_1\\
v_2
\end{pmatrix} \longleftarrow \begin{pmatrix}
A - \lambda_j I_n & B_2\\
C_1 & D_{12}
\end{pmatrix}\begin{pmatrix}
0_r\\
\hat{\eta}_j\\
\xi_j
\end{pmatrix}$\;
$a_j \longleftarrow v_1 + \lambda_j e_j$, $c_j \longleftarrow v_2$\;
\tcc{$\begin{pmatrix}
-e_j\\
\hat{\eta}_j\\
\xi_j
\end{pmatrix}\in\mathbb{R}^{n+m_2}$ is the null vector associated with zero $\lambda_j$. }
}
\Return{$A$, $C_1$, $B_1$, $B_2$, $D_{11}$ and $D_{12}$}\;
\caption{Algorithm to generate control system (\ref{system}) that has given  stable invariant zeros}
\end{algorithm}
\ULforem

The control system (\ref{system}) generated by Algorithm \ref{algo} has stable invariant zeros $\lambda_1, \ldots, \lambda_r$. In fact, we have 
\begin{align*}
A\begin{pmatrix}
-e_j \\
\hat{\eta}_j\\
\end{pmatrix} + B_2\xi_j &= -v_1 -\lambda_j e_j + \sum_{k=r+1}^n a_k \hat{\eta}_{jk} + B_2\xi_j\\
&= \lambda_j \hat{\eta}_j - \sum_{k=r+1}^n a_k\hat{\eta} _{jk} - B_2\xi_j -\lambda_j e_j + \sum_{k=r+1}^n a_k \hat{\eta}_{jk} + B_2\xi_j= \lambda_j\begin{pmatrix}
-e_j \\
\hat{\eta}_j\\
\end{pmatrix}\\
C_{1}\begin{pmatrix}
-e_j \\
\hat{\eta}_j\\
\end{pmatrix} + D_{12}\xi_j &= -v_2 +\sum_{k=r+1}^n c_k\hat{\eta}_{jk} +D_{12}\xi_j\\
&= -\sum_{k=r+1}^n c_k\hat{\eta}_{jk} -D_{12}\xi_j + \sum_{k=r+1}^n c_k\hat{\eta}_{jk} +D_{12}\xi_j = 0. 
\end{align*}
Here $e_j\in\mathbb{R}^r$ is the $j$th unit vector, and $a_k$ and $c_k$ are the $k$th column vectors of $A$  and $C_1$, respectively. This result implies that $(-e_j^T, \hat{\eta}_j^T, \xi_j^T)^T$ is the null vector associated with stable invariant zero $\lambda_j$. Therefore the control system generated  by Algorithm \ref{algo} has stable zeros $\lambda_1, \ldots, \lambda_r$. 

\begin{figure}[htbp]
\centering
\includegraphics[width=5cm]{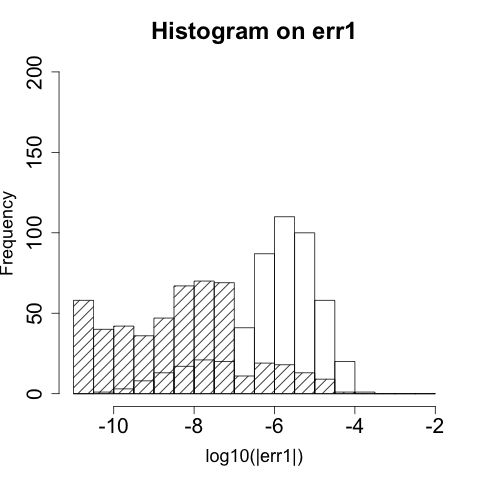}
\includegraphics[width=5cm]{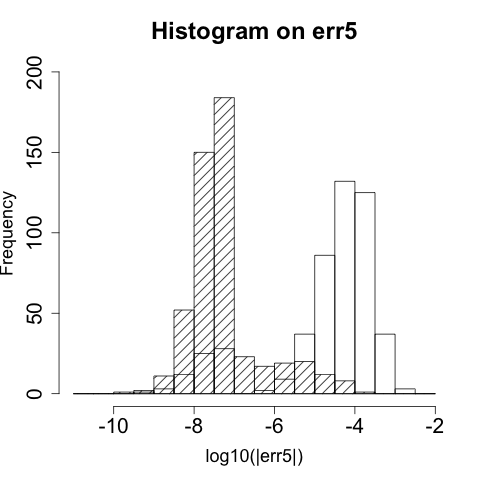}
\includegraphics[width=5cm]{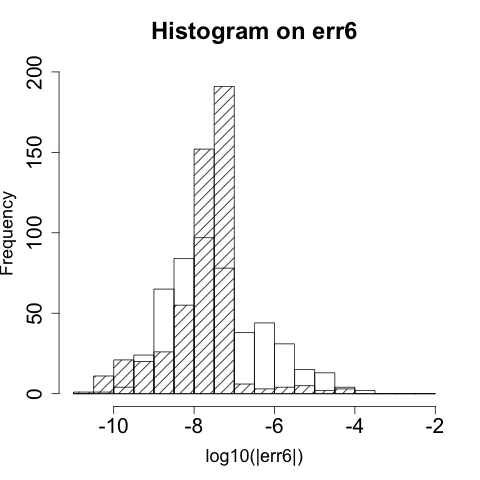}
\caption{Histograms of $\log_{10}(|\mbox{err}_1|)$ (left), $\log_{10}(|\mbox{err}_5|)$ (center) and $\log_{10}(|\mbox{err}_6|)$ (right) for (\ref{LMI2}) (plain boxes) and (\ref{LMI3}) (shaded boxes)  at $(n,  p_1, m_1, m_2) = (7, 5, 5, 2)$}
\label{fig6}
\includegraphics[width=5cm]{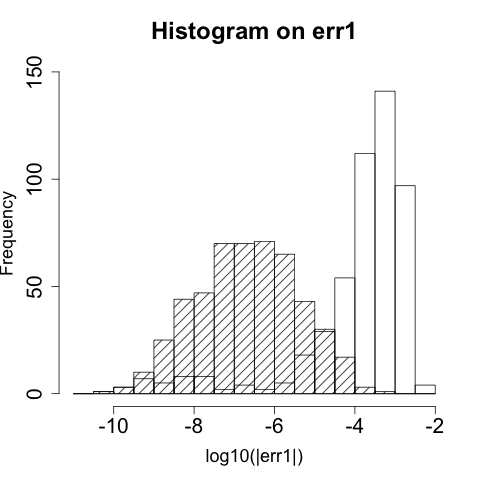}
\includegraphics[width=5cm]{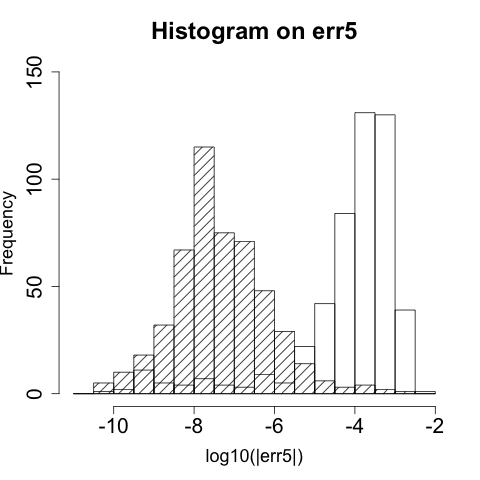}
\includegraphics[width=5cm]{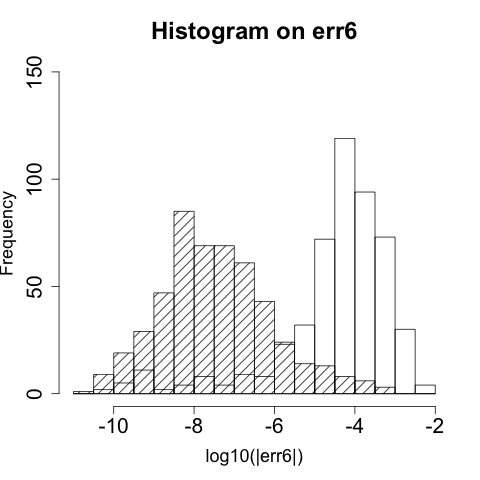}
\caption{Histograms of  $\log_{10}(|\mbox{err}_1|)$ (left), $\log_{10}(|\mbox{err}_5|)$ (center) and $\log_{10}(|\mbox{err}_6|)$ (right) for (\ref{LMI2}) (plain boxes) and  (\ref{LMI4}) (shaded boxes) at $(n,  p_1, m_1, m_2) = (7, 5, 5, 2)$}
\label{fig9}
\end{figure}

Figure \ref{fig6} displays histograms on the logarithms of the absolute values of $\mbox{err}_1$, $\mbox{err}_5$ and $\mbox{err}_6$ for the computed solutions of LMI problems (\ref{LMI2}) and (\ref{LMI3}). We generate 500 control  systems with three stable invariant zeros $\lambda_1 = -1$, $\lambda_2= -2$, $\lambda_3= -3$ and $(n,  p_1, m_1, m_2) = (7, 5, 5, 2)$ by Algorithm \ref{algo}. Figure \ref{fig9} displays histograms on  the logarithms of the absolute values of $\mbox{err}_1$, $\mbox{err}_5$ and $\mbox{err}_6$ for the computed solutions of LMI problems (\ref{LMI2}) and (\ref{LMI4}). We  randomly generate 500 control systems (\ref{system}) whose $D_{12}$ has the form $(O_{p_1\times 1}, \hat{D}_{12})$ for some $\hat{D}_{12}\in\mathbb{R}^{p_1\times (m_2-1)}$. The plain boxes in these figures indicate the results obtained by solving LMI obtained from (\ref{system}), while the shaded boxes indicate the results by (\ref{LMI3}) or (\ref{LMI4}). Table \ref{table:neg} displays the number of negative $\mbox{err}_5$ for the computed solutions. Since $\mbox{err}_5$ corresponds to the duality gap on the computed solution, the solution is not feasible if the value is negative.

\begin{table}[htbp]
\caption{The numbers of negative $\mbox{err}_5$}
\centering
\begin{tabular}{|c|c|c|c|c|}
\hline
& & $\#$ of $\mbox{err}_5 < -1.0\times 10^{-7}$ & $\#$ of $\mbox{err}_5 < -1.0\times 10^{-5}$ & $\#$ of $\mbox{err}_5 < -1.0\times 10^{-3 }$\\ 
\hline
Figure \ref{fig6} &(\ref{LMI2}) & 425 / 500&  381 / 500& 3 / 500 \\
&(\ref{LMI3}) & 77 / 500& 16 / 500& 0 / 500\\
\hline
Figure \ref{fig9} &(\ref{LMI2}) & 459 / 500& 423 / 500& 40 / 500\\
&(\ref{LMI4}) & 162 / 500 & 16 / 500& 1 / 500\\
\hline
\end{tabular}
\label{table:neg}
\end{table}%

 We observe the followings from figures:
\begin{itemize} 
\item Since the dual (\ref{dualLMI2}) of LMI (\ref{LMI2}) obtained from (\ref{system}) is not strongly feasible, the absolute values of the dual feasibility $\mbox{err}_1$ of (\ref{LMI2}) are worse than (\ref{clsubsystem}) and (\ref{rsss}). In particular, the computed dual solution of (\ref{LMI2}) is positive semidefinite, while it does not satisfy equality constraints in (\ref{LMI2}).   In fact, most of all the absolute values of the errors of (\ref{LMI3}) is smaller than the stopping tolerance $\epsilon = 1.0\times 10^{-7}$, while most of the absolute values of the errors $\mbox{err}_1$ and $\mbox{err}_5$ of (\ref{LMI2}) are bigger than the stopping tolerance. 
This means that we can obtain more accurate dual solutions by solving (\ref{LMI3}).  
\item The relative complementarity $\mbox{err}_5$ is the error on the duality gap.  Since the dual solution of (\ref{LMI2}) is not accurate, the computed dual solution of (\ref{LMI2}) has much worse feasibility than (\ref{LMI3}) and (\ref{LMI4}). Table \ref{table:neg} displays the numbers of negative $\mbox{err}_5$ of the computed solutions. In fact, if $\mbox{err}_5$ is negative, the weak duality fails in the computed solution, which implies that the computed dual solution is not feasible. In particular,  since the number of negative $\mbox{err}_5$ of (\ref{LMI2}) whose absolute value is the stopping tolerance is more than (\ref{LMI3}) and (\ref{LMI4}), the dual solution of (\ref{LMI2}) is less accurate than others. 

\item In the case where control system (\ref{LMI2}) has stable invariant zeros, the absolute values of  $\mbox{err}_6$ of (\ref{LMI2}) are similar to (\ref{LMI3}). One of the reasons may be that the computed primal solution of (\ref{LMI2}) is feasible and the zero eigenvalues of the dual solution are accurately computed. In contrast, the absolute values of $\mbox{err}_6$ are  bigger than the stopping tolerance $\epsilon$ in the case where $D_{12}$ is not full column rank. This implies that it is numerically difficult to solve LMI problems obtained from $H_{\infty}$ state feedback control problems in this case. 

\item Even if we reduce the control system via the stable invariant zeros, the errors for some of computed solutions of (\ref{LMI3})  and (\ref{LMI4}) are bigger than the stopping tolerance $\epsilon$. Moreover, $\mbox{err}_5$ is negative and its absolute value is bigger than $\epsilon$.  These mean that our proposed reduction is not sufficient to improve the numerical accuracy of $H_{\infty}$ state feedback control problem.  The improvement of the numerical accuracy is involved in future work. 
\end{itemize}

\section{conclusion}\label{sec:conclusion}

We discuss a numerical difficulty in solving LMI problem (\ref{LMI2}) and its dual (\ref{dualLMI2}) obtained from $H_{\infty}$ state feedback control, and show that the dual (\ref{dualLMI2}) is not strongly feasible if system (\ref{system}) has invariant zeros in $\overline{\mathbb{C}_-}$. This is derived from the viewpoint of facial reduction into LMI problems. Moreover, facial reduction provides the transformation with $T$ to reduce the size of system (\ref{clsystem}). 
 We observe in numerical results that the numerical stability in solving the resulting LMI problem  and its dual  is improved.

 This is not complete understanding of the numerical difficulty in solving LMI problems. In fact, DIMACS errors to  (\ref{LMI3}) and (\ref{LMI4})  sometimes become worse than the stopping tolerance $\epsilon$ of PDIPMs as in the table and figures in Section \ref{sec:experiment}. In addition, we often see the numerical difficulty in LMI problems obtained from state feedback control of system which does not have any invariant zeros. A more improvement of the numerical accuracy  for $H_{\infty}$ control problems are future work.

\section*{Acknowledgements}\label{sec:acknowledgements}
The first author was supported by JSPS KAKENHI Grant Numbers 22740056 and 26400203. We would like to thank  Dr. Yoshio Ebihara in Kyoto Univ. for a fruitful discussion and significant comments for improving the presentation of the manuscript. 

\appendix
\section{How to take the dual of  (\ref{LMI2})}\label{sec:dual}

The dual (\ref{dualLMI2}) of (\ref{LMI2}) can be obtained by reformulating (\ref{LMI2}) to (\ref{LMI}). We provide the detail in this section. 

Let $\{E_{ij}\}_{1\le i\le j\le n}$ and $\{F_{k\ell}\}_{1\le k\le m_2, 1\le \ell\le n}$ be base of the spaces $\mathbb{S}^{n}$ and $\mathbb{R}^{m_2\times n}$, respectively. We denote $X$ and $Y$ in (\ref{LMI2}) by $X = \sum_{1\le i\le j\le n}x_{ij}E_{ij}$ and $Y=\sum_{k=1}^{m_2}\sum_{\ell=1}^{n}y_{k\ell}F_{k\ell}$. Substitute them to (\ref{LMI2}), we can rewrite the constraints in (\ref{LMI2}) to 
\begin{align}\label{LMI2_tmp}
& \sum_{ij}x_{ij}\begin{pmatrix}
-\He(AE_{ij}) & -E_{ij}C_1^T & &\\
-C_1E_{ij} & & &\\
& & &\\
&&&E_{ij}
\end{pmatrix} + \sum_{k\ell}y_{k\ell}\begin{pmatrix}
-\He(B_2F_{k\ell}) & -F_{k\ell}^TD_{12}^T & &\\
-D_{12}F_{k\ell} & &&\\
& & & \\
& & & \\
\end{pmatrix} \nonumber\\
&{}+ \gamma \begin{pmatrix}
& & &\\
&I_{p_1}& &\\
& &I_{m_1} & \\
& & & \\
\end{pmatrix} - \begin{pmatrix}
& &B_{1} &\\
&&D_{11} &\\
B_{1}^T&D_{11}^T& & \\
& & & \\
\end{pmatrix}\in\mathbb{S}^{N_0+n}_+. 
\end{align}
Hence the dual is formulated as 
\begin{align}\label{dualLMI2_tmp}
&\left\{
\begin{array}{cl}
\sup & \begin{pmatrix}
& &B_{1} &\\
&&D_{11}&\\
B_{1}^T&D_{11}^T&&\\
& & & 
\end{pmatrix}\bullet \begin{pmatrix}
Z &V^T \\
V & W
\end{pmatrix}\\
\mbox{subject to} & W\bullet E_{ij} - \He(AE_{ij})\bullet Z_{11} - C_1E_{ij}\bullet Z_{21} - E_{ij}C_1^T\bullet Z_{21}^T = 0 \ (1\le i\le j\le n), \\
& -\He(B_2F_{k\ell})\bullet Z_{11} - D_{12}F_{k\ell}\bullet Z_{21} - F_{k\ell}^TD^T_{12}\bullet Z_{21}^T = 0 \ (1\le k\le m_2, 1\le\ell\le n), \\
& Z_{22}\bullet I_{p_1} + Z_{33}\bullet I_{m_1} = 1\\
& \begin{pmatrix}
Z &V^T \\
V & W
\end{pmatrix}\in\mathbb{S}^{N_0+n}_+\\
\end{array}
\right. 
\end{align}
Here $Z$ is partitioned as in (\ref{dualLMI2}). We can set $V=O_{n\times N_0}$ because $V$ does not appear in the objective function and equality constraints in (\ref{dualLMI2_tmp}). In addition, the equality constraints are equivalent to 
\begin{align*}
E_{ij}\bullet\left(
W - \He(A^TZ_{11} + C_1^TZ_{21}
\right) &= 0, F_{k\ell}\bullet \left(
\He(B_2^TZ_{11} +D_{12}^TZ_{21}
\right) = 0. 
\end{align*}
Since $E_{ij}$ and $F_{k\ell}$ are base of $\mathbb{S}^{n}$ and $\mathbb{R}^{m_2\times n}$, respectively, (\ref{dualLMI2_tmp}) is equivalent to (\ref{dualLMI2}). 

\section{Proofs}\label{sec:appendix}
\subsection{Proofs of the if-part and the infeasibility in Theorem \ref{nonsf}}\label{subsec:nonsf}
Suppose that there exists a nonzero $\hat{Y}\in\mathbb{S}^n_+$ such that $F_j\bullet \hat{Y} = 0 \ (j=1, \ldots, m)$ and $F_0\bullet \hat{Y} \ge 0$. If (\ref{LMI}) is strongly feasible, then there exists $(\tilde{x}, \tilde{X})\in\mathbb{R}^m\times\mathbb{S}^n_{++}$ such that $\tilde{X}=\displaystyle\sum_{j\in\mathcal{M}}\tilde{x}_jF_j-F_0$. Using $\hat{Y}$ and $(\tilde{x}, \tilde{X})$, we obtain
\[
0< \tilde{X}\bullet\hat{Y} = \left(\sum_{j\in\mathcal{M}}\tilde{x}_jF_j-F_0\right)\bullet \hat{Y} =  - F_0\bullet \hat{Y}\le 0. 
\]
Here, the first strict inequality is due to the positive definiteness of $\tilde{X}$. 
This implies the contradiction, and thus (\ref{LMI}) is not strongly feasible. 

We prove the infeasibility. Suppose  that $\hat{Y}\in\mathbb{S}^n_+$ satisfies $F_j\bullet\hat{Y}=0 \ (j\in\mathcal{M})$ and $F_0\bullet \hat{Y} > 0$. If (\ref{LMI}) is feasible, then we have $(\tilde{x}, \tilde{X})\in\mathbb{R}^m\times\mathbb{S}^n$ such that $\tilde{X}=\displaystyle\sum_{j\in\mathcal{M}}\tilde{x}_jF_j-F_0, \tilde{X}\in\mathbb{S}^n_{+}$. We obtain 
\[
0\le\tilde{X}\bullet\hat{Y} = \left(\sum_{j\in\mathcal{M}}\tilde{x}_jF_j-F_0\right)\bullet \hat{Y} =  - F_0\bullet \hat{Y}< 0, 
\]
and thus this implies the contradiction. Therefore, (\ref{LMI}) is infeasible. We can prove the remainder of this statement by applying a similar discussion to (\ref{dualLMI}). 

\subsection{Proof of Proposition \ref{propPNSF}}\label{sec:propPNSF}
We use the following lemma to prove some facts including Proposition \ref{propPNSF}: 

\lemm\label{techlemma}(See \cite[Lemma 2.4]{Ebihara2012})
For given $F, G\in\mathbb{M}^{n\times r}$, suppose that $F$ is full column rank. Then $FG^T+GF^T\in\mathbb{S}^n_+$ if and only if there exists $\Omega\in\mathbb{R}^{r\times r}$ such that $G=F\Omega$ and $\Omega + \Omega^T\in\mathbb{S}^{r}_+$. 
\elemm

We provide a proof of Proposition \ref{propPNSF}. To prove (only-if-part), suppose that problem (\ref{PNSFsystem2}) has a nonzero solution. Then, there exists a full column rank matrix $H\in\mathbb{R}^{n\times r}$ such that $Z_{11}=HH^T$, where $r$ is the rank of $Z_{11}$. (\ref{PNSFsystem2}) can be reformulated as follows: 
 \[
 A^TH H^T + HH^TA \in\mathbb{S}^{n}_+, B_2^TH = O_{m_2\times r}. 
 \] 
 It follows from the first inequality and Lemma \ref{techlemma} that there exists $\Omega\in\mathbb{R}^{r\times r}$ such that $\Omega + \Omega^T \in\mathbb{S}^r_+$ and $A^TH = H\Omega$.  This implies that all the eigenvalues of $\Omega$ are in $\overline{\mathbb{C}_+}$. Let $\eta$ be an eigenvector of $\Omega$. Then, we have $A(H\eta) = \lambda (H\eta)$, $\lambda\in\overline{\mathbb{C}_+}$ and $B_2\eta = 0$. Therefore, we obtain (\ref{stabrank}). 
 
 To prove (if-part), suppose (\ref{stabrank}) holds. Then there exist  $\lambda\in\overline{\mathbb{C}_+}$ and $\eta\in\mathbb{C}^{n}\setminus\{0\}$ such that $A^T\eta = \lambda \eta$ and $B_2^T\eta=0$. We define $Z_{11} = (\eta\bar{\eta}^T + \bar{\eta}\eta^T)$. Then $Z_{11}$ is positive semidefinite and we have 
 \begin{align*}
 A^TZ_{11} + Z_{11}A^T &= A^T\eta\bar{\eta}^T + A^T\bar{\eta}\eta^T + \eta\bar{\eta}^TA + \bar{\eta}\eta^TA\\
 &= (\lambda+ \bar{\lambda})(\eta\bar{\eta}^T + \bar{\eta}\eta^T) = (\lambda+ \bar{\lambda})Z_{11}\in\mathbb{S}^n_+, \\
 B_2^TZ_{11} &= B_{2}(\eta\bar{\eta}^T + \bar{\eta}\eta^T)=B_{2}\bar{\eta}\eta^T = 0. 
 \end{align*}
 The last equality holds because we have $B_{2}\bar{\eta}= 0$ by taking the complex conjugate to $B_2\eta = 0$. These imply that  (\ref{PNSFsystem2}) holds.

\subsection{Proof of Proposition \ref{propDNSF}}\label{sec:propDNSF}
We provide a proof of Proposition \ref{propDNSF}. Suppose that there exists $\lambda\in\overline{\mathbb{C}_-}$ which satisfies (\ref{rankinvz}). 
Then, we have $\eta\in\mathbb{C}^n\setminus\{0\}$ and $\xi\in\mathbb{C}^{m_2}$ satisfying $A\eta + B_{2}\xi = \lambda \eta, C_1\eta + D_{12}\xi = 0$. We define $X = (\eta\bar{\eta}^T + \bar{\eta}\eta^T)$ and $Y = (\xi\bar{\eta}^T + \bar{\xi}\eta^T)$. Then, $X$ is positive semidefinite and we have 
\begin{align*}
 C_1X + D_{12} Y &= (C_1\eta+D_{12}\xi)\bar{\eta}^T + (C_1\bar{\eta} + D_{12}\bar{\xi})\eta^T = 0\\
 AX+B_2Y &= (A\eta +B_2\xi)\bar{\eta}^T + (A\bar{\eta} + B_2\bar{\xi})\eta^T = \lambda \eta\bar{\eta}^T +\bar{\lambda}\bar{\eta}\eta^T \\
 -\He(AX+B_2Y) &=-(\lambda+\bar{\lambda})(\eta\bar{\eta}^T +\bar{\eta}\eta^T) = -(\lambda+\bar{\lambda})X\in\mathbb{S}^n_+,   
\end{align*}
which implies that (\ref{DNSFsystem2}) holds. 

Suppose that there exists a nonzero solution $(X, Y)$ of (\ref{DNSFsystem2}) such that $X=HH^T$ and $Y=RH^T$ for a full column matrix $H\in\mathbb{R}^{n\times r}$ and $R\in\mathbb{R}^{p_1\times r}$. Then we have 
\[
C_1H +D_{12}R = O_{p_1\times r}, \He((-AH - B_2R)H^T) \in\mathbb{S}^n_+. 
\]
It follows from the last inequality and Lemma \ref{techlemma} that there exists $\Omega\in\mathbb{R}^{r\times r}$ such that $-(AH +B_2R) = H\Omega$  and  $\Omega + \Omega^T\in\mathbb{S}^{r}_+$.
All eigenvalues of $\Omega$ are in $\overline{\mathbb{C}_+}$. Let $\eta$ be an eigenvector corresponding to an eigenvalue $\lambda$ of $\Omega$. Then we have
$C_1(-H\eta) +D_{12}(-R\eta) = 0$ and
 $(-H\eta) +B_2(-R\eta) = -\lambda (-H\eta)$. 
$-H\eta$ is nonzero because $H$ is full column rank, which implies that (\ref{rankinvz}) holds. 

\section{Computation for dual (\ref{dualLMI2_fr2})}\label{appendix:computation}
We give a detail how to obtain (\ref{dualLMI2_fr2}). To this end, we need to compute $\tilde{A}, \tilde{B}_1, \tilde{B}_2$ and $\tilde{C}_1$ in (\ref{dualLMI2_fr2}). 

 For the constraint $\He(A^TZ_{11}+C_{1}^TZ_{21})\in\mathbb{S}^n_+$ in  (\ref{dualLMI2_fr}), multiplying $T$ and $T^T$, then we obtain $\He((T^{-1}AT)^T(T^TZ_{11}T) +(C_{1}T)^TZ_{21}T))\in\mathbb{S}^n_+$.  In addition, we have 
 \begin{align*}
T^{-1}AT &= \begin{pmatrix}
\tilde{A}_{11} & \tilde{A}_{12}\\
\tilde{A}_{21} & \tilde{A}_{22}
\end{pmatrix}, C_1T = \begin{pmatrix}\tilde{C}_{11}, \tilde{C}_{12}\end{pmatrix}, \\
 T^T Z_{11}T &= \begin{pmatrix}
O_{r\times r} & O_{r\times (n-r)}\\
O_{(n-r)\times r}&\tilde{Z}_{11}\end{pmatrix}, Z_{21}T=\begin{pmatrix}
O_{p_1\times r} & \tilde{Z}_{21}
\end{pmatrix}.
 \end{align*}
 By a direct computation, the constraint is equivalent to 
 \[
 \He\left(\begin{pmatrix}
 O_{r\times r} & \tilde{A}_{21}^T\tilde{Z}_{11} +\tilde{C}_{11}^T\tilde{Z}_{21}\\
 O_{(n-r)\times r} & \tilde{A}_{22}^T\tilde{Z}_{11} + \tilde{C}_{12}^T\tilde{Z}_{21}
 \end{pmatrix}\right)\in\mathbb{S}^n_+. 
 \] 
 From this constraint, we obtain $\tilde{A}_{21}^T\tilde{Z}_{11} +\tilde{C}_{11}^T\tilde{Z}_{21} = O_{r\times (n-r)}$ and $\He(\tilde{A}_{22}^T\tilde{Z}_{11} + \tilde{C}_{12}^T\tilde{Z}_{21})
 \in\mathbb{S}^{n-r}_+$.  
 
 For the constraint $B_2^TZ_{11} + D_{12}^TZ_{21} = O_{m_2\times n}$, multiplying $T$, we obtain $B_2^TT^{-T}T^TZ_{11}T+ D_{12}^TZ_{21}T = O_{m_2\times n}$. Then, we have 
\[
 \begin{pmatrix}
 \tilde{B}_{21}^T &  \tilde{B}_{22}^T 
 \end{pmatrix}\begin{pmatrix}
O_{r\times r} & O_{r\times (n-r)}\\
O_{(n-r)\times r} & \tilde{Z}_{11}
\end{pmatrix} + D_{12}^T \begin{pmatrix}
O_{p_1\times r} & \tilde{Z}_{21} 
\end{pmatrix} = O
 \]
 and thus $\tilde{B}_{22}^T \tilde{Z}_{11}+ D_{12}^T\tilde{Z}_{21} = O_{m_2\times (n-r)}$. 
 We can obtain  dual problem (\ref{dualLMI2_fr2}) combining these computation with (\ref{reductionL}). 
\end{document}